# Categorical Aspects in Pasch Geometry

Hidayet Hüda KÖSAL[†] and Emine SOYTÜRK SEYRANTEPE[‡]

**Abstract.** In this study, after introducing Pasch geometries and algebraic properties of their, we studied categorical structure of the categories of Pasch geometries with morphisms and homomorphisms. In this regard, definition and construction of some categorical objects such as: equalisers, pullbacks, initial and terminal objects and products are given.

1. **Introduction**

The basic theory of Pasch geometries has been developed by Harrison, D.K. [1]. Also Harrison, D.K. defined morphism and homomorphisms between Pasch geometries. These morphisms and homomorphisms can be a partially defined map between the points sets, satisfying simple geometric conditions and form categories.

The geometries of double cosets and orbits of groups are studied abstractly as Pasch geometries, hypergroups ([2], [3]), multigroups ([4]). In particular, the points of the projective space $P(V)$ of a vector space $V$ over a field $F$ can be taken as the set of orbits $V/F^*$ and the inherited structure provides the geometry of the projective space. In addition, some of these structures posses probabilistic measures whereby there is a certain probability for a given element to belong to the set of products of two given elements. Such structures have been abstractly studied as Probability groups ([5], [1]). It has been proved in [6] that an elementary abelian Pasch geometry is isomorphic as geometries to an orbit space of an abelian group.

Morphisms and Homomorphisms of Projective geometries are introduced in [7]. Categories of Projective geometries with morphisms and homomorphisms has been studied by Bhattarai in [8]. Projective geometries form categories with morphism as well as homomorphism and Desarguesian ones form a subgeometry with Desarguesian homomorphism. In addition categories of projective geometries are full subgeometry of categories of Pasch geometries. This paper will focus on definition and contruction of some categorical objects of categories of Pasch geometries.

The paper is organized in the following way: in sect. 2 the algebraic properties of Pasch geometries are summarized. In sect. 3 some categorical objects of the categories of Pasch geometries are given.

[†] Sakarya University, Faculty of Arts and Sciences Department of Mathematics, Sakarya / Turkey. E-mail: hhkosal@sakarya.edu.tr
[‡] Afyon Kocatepe University, Faculty of Arts and Sciences Department of Mathematics, Afyon/ Turkey. E-mail: soyturk@aku.edu.tr

## 2. Algebraic Properties of Pasch Geometries

The basic concepts of Pasch geometry can be references, particulary in [1] and [2]. For convenience we briefly state some relevant definitions and result here.

**Definition 2.1.** By a Pasch geometry it is meant a triple $(A, e, \Delta_A)$ where $A$ is a set, $e$ is an element in $A$, and $\Delta_A$ is a subset of $A \times A \times A$ subject to the following axioms:

(1) For each $a \in A$ there exist a unique $b \in A$ with $(a, b, e) \in \Delta_A$; denote $b$ by $a^\#$.

(2) $e^\# = e$ and $(a^\#)^\# = a$ for all $a \in A$.

(3) $(a, b, c) \in \Delta_A$ implies $(b, c, a) \in \Delta_A$.

(4) $(a_1, a_2, a_3), (a_1, a_4, a_5) \in \Delta_A$ imply there exists an $a_6 \in A$ with $(a_6, a_4^\#, a_2), (a_6, a_5, a_3^\#)$.

As a consequence of the above one gets:

(5) $(a, b, c) \in \Delta_A$ implies $(c^\#, b^\#, a^\#) \in \Delta_A$.

(6) $x, y \in A$ implies $\exists z \in A$ with $(x, y, z) \in \Delta_A$.

Throughout this study, we will sometimes simply write "geometry" for "Pasch Geometry".

A geometry is called Abelian if $(a, b, c) \in \Delta_A$ implies $(b, a, c) \in \Delta_A$. We say a Pasch geometry $(A, e, \Delta_A)$ is sharp if for each $a, b \in A$ there exists at most one $c \in A$ with $(a, b, c) \in \Delta_A$. One can show there always is at least one such $c$, and if $a.b$ denotes $c^\#$, one checks a group result. Conversely every group $G$ is naturally a sharp geometry with $(x, y, z) \in \Delta_G$ if and only if $xyz = e$, the identity of $G$.

Let $A$ be a geometry and $S \subseteq A$. Then $S$ is called a subgeometry if $e \in S$ and $s_1, s_2 \in S$, $(s_1, s_2, x) \in \Delta_A$ implies $x \in S$. A subgeometry $S$ is called normal if $\forall a, b \in A$, $(s, a, b) \in \Delta_A$ for some $s \in S$ implies $\exists s_1 \in S$ with $(s_1, b, a) \in \Delta_A$.

Let $A$ and $B$ be geometries. A map $f : A \to B$ is called a geometry morphism if $f(e_A) = e_B$ and $(x, y, z) \in \Delta_A$ implies $(f(x), f(y), f(z)) \in \Delta_B$. For $f$ such a morphism let $K_f = \{a \in A : f(a) = e_B\}$, called the kernel of $f$ and $\operatorname{Im} f = \{b \in B : b = f(x), \text{ for some } x \in A\}$ called the image of $f$. Then $K_f$ is a subgeometry of $A$ which may not be normal. Also $\operatorname{Im} f$ may not be a subgeometry of the codomain.

A map $f : A \to B$ is called a homomorphism if it is a morphism and $(f(x), f(y), b) \in \Delta_B$ implies $\exists z \in A$ with $b = f(z)$ and $(x, y, z) \in \Delta_A$.

**Proposition 2.2.** Let $C$ is a geometry, $g$ is a morphism from $B$ to $C$, and $f$ is a morphism from $A$ to $B$. Then $g \circ f$ is a morphism from $A$ to $C$ [1].

**Proposition 2.3.** Let $A$ and $B$ be sharp geometries. If $f: A \to B$ is a geometry homomorphism, then $f$ is a group homomorphism [1].

**Proposition 2.4.** Let $A$ and $B$ be geometries. Then $(A \times B, (e_A, e_B), \Delta_{A \times B})$ is a geometry [1]. Where $\Delta_{A \times B}$ denotes

$$\{((a_1, b_1), (a_2, b_2), (a_3, b_3)) : (a_1, a_2, a_3) \in \Delta_A, (b_1, b_2, b_3) \in \Delta_B\}.$$

### 3. Categorical Structure of Pasch Geometries

The fundamental concept of Category Theory can be found in references, in particulary: Blyth, T.S., Categoires, 1986 [9]. Herlich, H. and Strecker, G.E., Category Theory, An Introduction, 1973 [10]. Higgins, P.,J., Notes on categories and groupoids, 1971 [11]. Barry M., Theory of Categories, 1965 [12].

The category $P$ : The objects of the category are Pasch Geometries. Morphisms *Mor(P)* is the set of all geometry homomorphisms of geometries. The set of all morphisms from $A$ to $B$ is denoted by $Mor(A, B)$.

**Theorem 3.1.** The category $P_1$ of Abelian geometries and all geometry homomorphisms between them is a subcategory of $P$.

**Theorem 3.2.** The category $P_2$ of sharp geometries and all geometry homomrphisms between them is a category of $P$.

**Theorem 3.4.** In the category $P$, each one-element geometry is both initial and terminal objects.

**Theorem 3.5.** Terminal object in the category $P$ is unique up to isomorphism.

**Proof.** Suppose that $T$ and $T'$ are terminal objects. Since $T$ is terminal, there is a unique morphism $f: T' \to T$. Similary, there is a unique morphism $g: T \to T'$. The morphism $f \circ g : T \to T$ is a morphism with target $T$. Since $T$ is a terminal object of the category, there can be only one morphism from $T$ to $T$. Thus it must be that $f \circ g$ is the identity of $T$. An analogous proof shows that $g \circ f$ is the identity of $T'$.

**Theorem 3.6.** Initial object in the category $P$ is unique up to isomorphism.

The proof is a similar check and is omitted.

**Corollary 3.7.** Since the one element geometry is both initial and terminal objects in the category of pasch geometries, It is a zero object in the category $P$.

**Theorem 3.6.** Let $A$ and $B$ two objects in the category $P$. Then product object of $A$ and $B$ is $A \times B$ where $A \times B$ denotes by $(A \times B, (e_A, e_B), \Delta_{A \times B})$.

**Proof.** $(A \times B, (e_A, e_B), \Delta_{A \times B})$ is a geometry (Proposition 2.4). Let $D$ be any geometry and $q_1 : D \to A$, $q_2 : D \to B$ be morphisms. Then, a map

$$q : D \to A \times B$$
$$x \to q(x) = (q_1(x), q_2(x))$$

is morphism. Firstly, we define

$$\pi_1 : A \times B \to A \qquad\qquad \pi_2 : A \times B \to A$$
$$(x, y) \to \pi_1(x, y) = x \qquad\qquad (x, y) \to \pi_2(x, y) = y$$

morphisms.

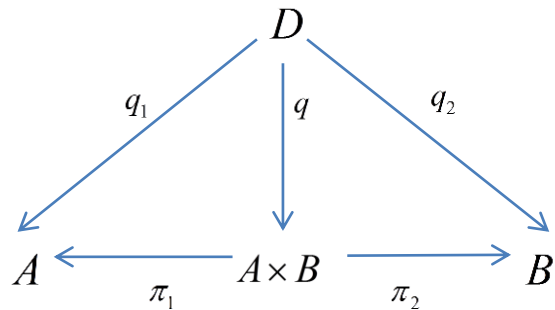

**Diyagram (1)**

Then, diyagram (1) is commutative. Now, we show that there is a unique $q : D \to A \times B$ morphism.

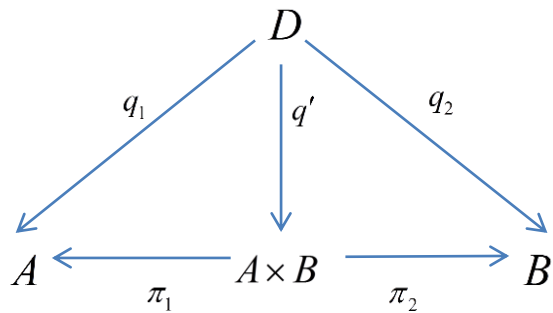

**Diyagram (2)**

Assume that,

$$q': D \to A \times B$$
$$x \to q'(x) = (q'_1(x), q'_2(x))$$

morphism make commutative to diyagram (2). Then, for $\forall x \in D$ $(\pi_1 \circ q')(x) = q_1(x)$ and $(\pi_2 \circ q')(x) = q_2(x)$ is satisfied. Due to

$$\pi_1(q'(x)) = \pi_1(q'_1(x), q'_2(x)) \qquad \pi_2(q'(x)) = \pi_2(q'_1(x), q'_2(x))$$
$$= q'_1(x) = q_1(x) \qquad\qquad = q'_2(x) = q_2(x),$$

we obtain $q = q'$. Thus, there exist a unique $q: D \to A \times B$ morphism.

**Theorem 3.7.** In the category $P_2$, let $f, g: A \to B$ be parallel pair of morphisms. Then an equalizer of $f$ and $g$ is an object $E = \{x | x \in A, f(x) = g(x)\}$ together with an morphism $i: E \to A$, where $i(x) = x$ for all elements $x \in E$.

**Proof.** Let $A, B$ be objects and $f, g: A \to B$ morphisms, in category $P$. Then, $E = \{x | x \in A, f(x) = g(x)\}$ is subgeometry of $A$. We define a function: $i: E \to A$, $i(x) = x$ for all elements $x \in E$. This function is morphism in the category $P_2$. Moreover, $f \circ i = g \circ i$ is satisfied. On the other hand, let $k: E' \to A$ be morphism for any $E'$ object. We assume that $f \circ k = g \circ k$. Due to $(f \circ k)(x) = f(k(x)) = g(k(x)) = (g \circ k)(x)$, we obtain $k(x) \in E$. We deduce that $\operatorname{Im}(k) \subset E$. Now, let taking $h: E' \to E$, $h(x) = k(x)$ for $\forall x \in E'$. $h$ is morphism and $i \circ h = k$ is satisfied. Finally, we deduce that there is a unique $h: E' \to E$. Assume that $h': E' \to E$ morphism make commutative to diyagram (3).

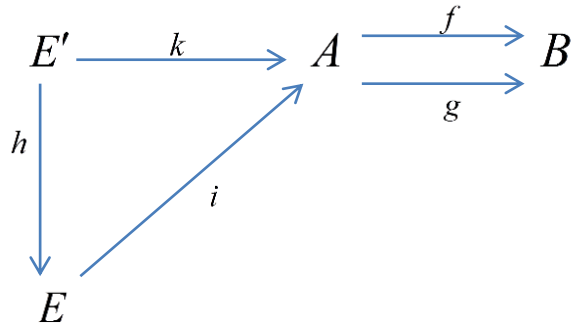

**Diyagram (3)**

Due to, $(h' \circ i)(x) = h'(i(x)) = h'(x) = k(x)$, for $\forall x \in E'$, we obtain $h = h'$. Thus, there exist a unique $h: E' \to E$.

Finally, equalizer of $f$ and $g$ is an object $E$ together with an morphism $i: E \to A$.

**Theorem 3.8.** In the category $P_2$, if $f: A \to X$ and $g: B \to X$ are morphism, then the pullback

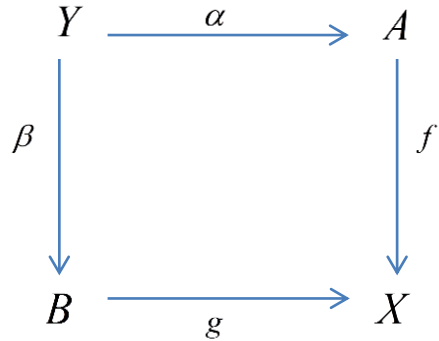

**Diyagram (4)**

is constructed by $Y = \{(a,b): f(a) = g(b)\}$ with $\alpha(a,b) = a$ and $\beta(a,b) = b$.

**Proof.** Let $f: A \to X$, $g: B \to X$ be morphisms and $Y = \{(a,b): f(a) = g(b)\}$. Then, $Y$ is subgeometry of $A \times B$. We define morphisms as follows:

$$\alpha: Y \to A \qquad \beta: Y \to B$$
$$(a,b) \to a \qquad (a,b) \to b.$$

It is clear that diyagram (4) is commutative. Now we show that diyagram (4) is a pullback diyagram.

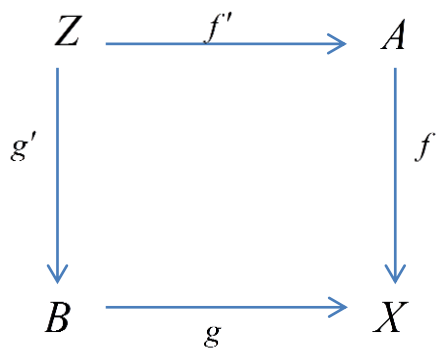

**Diyagram (5)**

We assume that diyagram (5) is commutative for any $Z$ object and $f': Z \to A$, $g': Z \to B$ morphisms. By using $(f \circ f')(x) = f(f'(x)) = g(g'(x)) = (g \circ g')(x)$, $\forall x \in Z$, we deduce that $(f'(x), g'(x)) \in Y$. We define a morphism

$$\varepsilon: Z \to Y$$
$$x \to \varepsilon(x) = (f'(x), g'(x)).$$

In this case, diyagram (6) is commutative. Moreover, there exist a unique $\varepsilon: Z \to Y$ morphism.

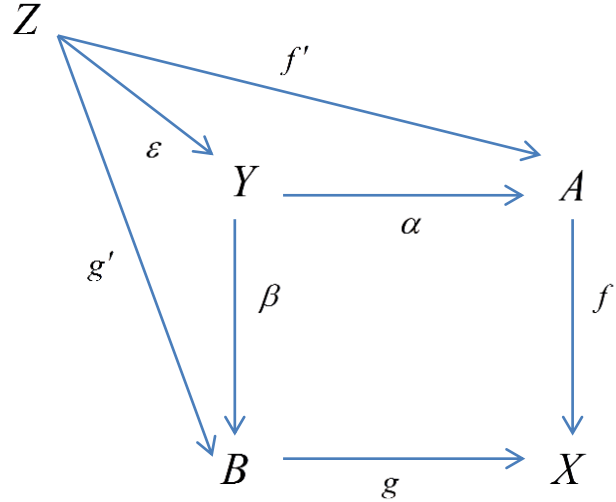

**Diyagram (6)**

**Theorem 3.9.** In the category $P$, $f \sim g$ is defined by

"$f \sim g$ if and only if there exist $b \in B$ such that $f(a) = bg(a)b^{\#}$ for all $a \in A$ and $f, g \in Mor(A, B)$."

This defines an equivalence relation on $Mor(P)$. $[f] = \{g | f \sim g\}$ for $f \in Mor(P)$ is called the set of equivalence class of $f$. Moreover, equivalence relation on $Mor(P)$ is congruence relation.

**Proof.** Firstly, we show that $\sim$ is an equivalence relation.

(i) $f(a) = f(a) f(a) [f(a)]^{\#}$ for any $a \in A$, thus $f \sim f$ and $\sim$ is reflexive.

(ii) Let $f \sim g$. This means there exist $b \in B$ such that $f(a) = bg(a)b^{\#}$ for $a \in A$. Let $d = b^{\#}$. Then, $df(a)d^{\#} = g(a)$. Hence $g \sim f$. So $\sim$ is symmetric.

(iii) Let $f \sim g$ and $g \sim h$. Then, there exist $b_1, b_2 \in B$ such that $f(a) = b_1 g(a) b_1^{\#}$ for all $\forall a \in A$ and $g(a) = b_2 h(a) b_2^{\#}$. Then

$$f(a) = b_1 g(a) b_1^{\#} = b_1 (b_2 h(a) b_2^{\#}) b_1^{\#} = (b_1 b_2) h(a) (b_2^{\#} b_1^{\#}) = (b_1 b_2) h(a) (b_1 b_2)^{\#}.$$

Thus $f \sim h$ and $\sim$ is transitive.

Now, we show that $\sim$ is a congruence relation.

(1) Let $f : A \to B$ is morphism and $[f]$ is an equivalence class of $f$. Then

$$[f] = \{g \mid g : A \to B, f \sim g\} \subseteq Mor(A, B).$$

(2) Let $f \sim f'$, $g \sim g'$ and $f, f' : A \to B$, $g, g' : B \to C$ be morphisms. Due to $f \sim f'$, there exist $b_1 \in B$ such that $f(a) = b_1 f'(a) b_1^{\#}$ for all $a \in A$. Similary, due to $g \sim g'$, there exist $c_1 \in C$ such that $g(b) = c_1 g'(b) c_1^{\#}$ for all $\forall b \in B$. Then

$$(g \circ f)(a) = g(f(a)) = g(b_1 f'(a) b_1^{\#}) = g(b_1) g(f'(a)) [g(b_1)]^{\#}$$
$$= g(b_1) c_1 g'(f'(a)) c_1^{\#} [g(b_1)]^{\#}$$
$$= (g(b_1) c_1)(g' \circ f')(a)(g(b_1) c_1)^{\#}, \text{ for all } a \in A.$$

Thus, $\sim$ equivalence is a congruence relation on $Mor(P)$.

Now, we construct the quatien category.

- $Ob(P/_\sim) = Ob(P)$

- $Mor(P/_\sim) = \{[f] \mid [f] \subseteq Mor(A, B)\}$

- $[f] \circ_{P/_\sim} [g] = [g \circ_P f]$, for $\bullet [f], [g] \in Mor(P/_\sim)$.

**Acknowledgement**